\documentclass[12pt]{article}
\usepackage[letterpaper,margin=1in]{geometry}
\usepackage{fancyhdr,color}
\usepackage{tikz,graphicx,multicol}
\usepackage{amssymb,euscript,nicefrac,enumitem}
\usepackage{amsfonts,amsmath,amsthm}
\usepackage{scrextend}
\usepackage{ytableau} 
    \ytableausetup{boxsize=1.0em,aligntableaux=center}
\usepackage{shuffle}
\usepackage{extarrows}
\usepackage{pdfpages}
\usepackage{mathrsfs}
\usepackage{arydshln}
\usepackage{algorithm}
\usepackage{algpseudocode}
\usetikzlibrary{cd}

\newtheorem{thm}{Theorem}

\newtheorem{defn}{Definition}
\newtheorem{example}{Example}

\usepackage[backend=bibtex,maxbibnames=99]{biblatex}
 \def\gdot{\color{black!25}{\scriptstyle\cdot}}
\usepackage{listings}

\title{Powersum Bases in Quasisymmetric Functions and Quasisymmetric Functions in Non-commuting Variables}
\author{Anthony Lazzeroni}
\date{}

\makeatletter
\let\theauthors\@author
\makeatother
\fancypagestyle{plain}{
  
  \fancyhf{}
  \fancyhead[L]{}
  \fancyhead[C]{}
  \fancyhead[R]{}
  \fancyfoot[C]{\footnotesize \thepage}
}

\def\sym{\mathrm{Sym}}
\def\qsym{\mathrm{QSym}}
\def\nsym{\mathrm{NSym}}
\def\ncsym{\mathrm{NCSym}}
\def\ncqsym{\mathrm{NCQSym}}
\def\fqsym{\mathrm{FQSym}}
\def\mqsym{\mathrm{MQSym}}

\def\des{\mathsf{des}}

\def\sort{\mathsf{sort}}
\def\sdr{\mathsf{SDR}}
\def\hgt{\mathsf{ht}}
\def\len{\mathsf{len}}
\def\row{\mathsf{row}}
\def\col{\mathsf{col}}

\def\qshuf{\underline{\shuffle}}
\def\sd{\mathsf{SD}}

\def\ldd{\mathsf{LDD}}

\begin{document}

\maketitle 

\abstract{We introduce new bases for the Hopf algebra of quasisymmetric functions that refine the symmetric powersum basis. These bases are expanded in terms of quasisymmetric monomial functions by using fillings of matrices. We define the analog of these bases in quasisymmetric functions of non-commuting variables. Our new bases have a (shifted) shuffle product and a deconcatenate coproduct. Finally, we describe a change of basis rule from the quasisymmetric powersum basis to the quasisymmetric fundamental basis.}

\section{Introduction}
The Hopf algebra of symmetric functions, denoted $\sym$ and indexed by partitions $\lambda$ and $\mu$, is a very well known space for its connections in representation theory and other areas of mathematics. Some of the well studied bases are the monomial basis $m_\lambda$, powersum basis $p_\lambda$, and Schur basis $s_\lambda$. In \cite{zabrocki2015introduction} Zabrocki illustrates the change of basis from $p_\lambda$ to $m_\mu$ by fillings. In this work, a \textit{filling} $\mathsf{F}$ is a matrix with one non-zero entry in every row. The column (respectively, row) reading is a composition recording the sum of all entries in each column (respectively, row) denoted as $\col(\mathsf{F})$ (respectively, $\row(\mathsf{F})$). Thus the change of basis is combinatorially defined as 
\begin{equation} \label{eqn:ptom}
    p_\lambda = \sum_{\mathsf{F}\in A(\lambda)}m_{\col(\mathsf{F})}
\end{equation}
where $A(\lambda)$ is the set of all the distinct fillings with row reading $\lambda$ and the column reading is a partition. For example, $p_{(2,2,1)}=2m_{(2,2,1)}+2m_{(3,2)}+m_{(4,1)}+m_{(5)}$ since all the fillings of $A(2,2,1)$ are
\[ \footnotesize
\begin{array}{c|ccc}
2 & 2 & \gdot  & \gdot  \\
2 & \gdot  & 2 & \gdot  \\
1 & \gdot  & \gdot  & 1 \\
\hline
& 2 & 2 & 1
\end{array}
\qquad
\begin{array}{c|ccc}
2 & \gdot  & 2 & \gdot  \\
2 & 2 & \gdot  & \gdot  \\
1 & \gdot  & \gdot  & 1 \\
\hline
& 2 & 2 & 1
\end{array}
\qquad
\begin{array}{c|ccc}
2 & 2      & \gdot  & \gdot  \\
2 & \gdot  & 2      & \gdot  \\
1 & 1      & \gdot  & \gdot  \\
\hline
& 3 & 2 & 
\end{array}
\qquad
\begin{array}{c|ccc}
2 & \gdot  & 2      & \gdot  \\
2 & 2      & \gdot  & \gdot  \\
1 & 1      & \gdot  & \gdot  \\
\hline
& 3 & 2 & 
\end{array}
\qquad
\begin{array}{c|ccc}
2 & 2     & \gdot  & \gdot  \\
2 & 2     & \gdot  & \gdot  \\
1 & \gdot & 1      & \gdot  \\
\hline
& 4 & 1 & 
\end{array}
\qquad
\begin{array}{c|ccc}
2 & 2  & \gdot  & \gdot  \\
2 & 2  & \gdot  & \gdot  \\
1 & 1  & \gdot  & \gdot  \\
\hline
& 5 &  & 
\end{array}. 
\]
One of the important properties of the powersum basis is the Murnaghan-Nakayama rule which illustrates the product rule of $s_\lambda$ and $p_\mu$ expanded in terms of Schur functions.

A space that contains $\sym$ is the Hopf algebra of quasisymmetric functions, $\qsym$, and is indexed by compositions $\alpha$ and $\beta$. This space was defined in \cite{IraM.Gessel1984MultipartiteFunctions} by using $P$-partitions, which generates the fundamental quasisymmetric functions $F_\alpha$. In \cite{Ballantine2020OnSums} the authors introduced two quasisymmetric powersum bases (i.e. a basis of $\qsym$ that refines $p_\lambda$) $\Phi$ and $\Psi$ whose duals are the $\mathbf{\Phi}$ and $\mathbf{\Psi}$ bases in $\nsym$ which was introduced in \cite{Gelfand1995NoncommutativeFunctions}. These two bases aren't defined combinatorially, but by using series. The $\Psi$ basis is most notable for the change of basis to the fundamental basis by using $P$-partitions as defined in \cite{Alexandersson2021P-PartitionsP-Positivity}. In \cite{aliniaeifard2021peak} the authors introduced the Shuffle basis $S_\alpha$, which also refines $p_\lambda$, is notable because $S_\alpha$ is an eigenvector under the theta map $\Theta$.

In this extended abstract we define a quasisymmetric powersum basis $P_\alpha$ combinatorially by using fillings in an analog of Equation \ref{eqn:ptom}. Alternatively, $P_\alpha$ can be defined using a subposet of the refinement poset $\mathcal{P}$ on compositions. Both the fillings and the subposet have generalizations that depend on a total order $\preceq$ on the parts of the compositions. Hence we can define a whole family of quasisymmetric powersum bases, denoted as $P^{\succeq}_\alpha$ for different choices of $\succeq$, such that every basis has a shuffle product, a deconcatenate coproduct, and they refine the symmetric powersum basis. When $\succeq$ is the usual order $\geq$, our quasisymmetric powersum basis is dual (up to scaling) to the Zassenhaus functions $Z_\alpha$ in $\nsym$, the algebra of non-commutative symmetric functions, as defined in \cite{krob1997noncommutative}. All these results are stated in Section \ref{sec:qsym}.
We note that the special cases $P^{\geq}_\alpha$ and $P^{\leq}_\alpha$ were independently defined in \cite{aliniaeifard2021p} using weighted generating functions of P-partitions.

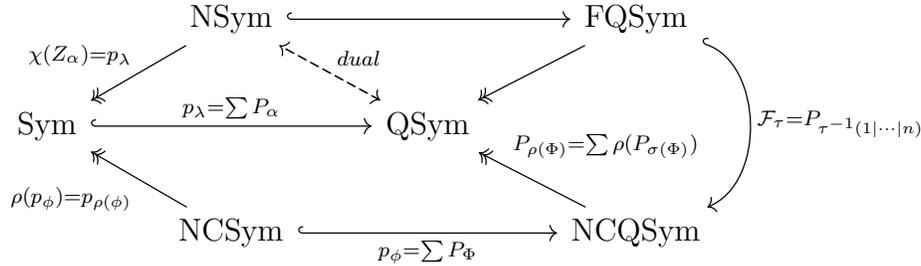
\begin{figure}[h] 
\begin{center}
\begin{tikzcd} 
 & \nsym \arrow[rr, hook] \arrow[dl, twoheadrightarrow, "\chi(Z_\alpha)=p_\lambda"'] & & \fqsym  \arrow[dl, twoheadrightarrow] \arrow[dd, hook, bend left=70, "\mathcal{F}_\tau=P_{\tau^{-1}(1|\cdots|n)}" ] \\
\sym \arrow[rr, hook, "p_\lambda = \sum P_\alpha"] & & \qsym \arrow[ul, dashleftarrow] \arrow[ul, dashrightarrow, "dual"'] \arrow[dr, twoheadleftarrow, "P_{\rho(\Phi)} = \sum \rho(P_{\sigma(\Phi)})" near start ]\\
 & \ncsym \arrow[ul, twoheadrightarrow, "{\rho(p_\phi)=p_{\rho(\phi)}}"] \arrow[rr, hook, "p_\phi = \sum P_\Phi"'] & & \ncqsym
\end{tikzcd}
\end{center}
\caption{Diagram of some Hopf Algebras related to $\qsym$}
\label{fig:alg}
\end{figure}

In Section \ref{sec:ncqsym}, we extend the above results to non-commuting variables - that is, we define a family of powersum bases $P_\Phi$ in $\ncqsym$, which refine the symmetric powersum basis $p_\phi$ of $\ncsym$ as defined in \cite{Rosas2006SymmetricVariables}. These bases have a shifted shuffle product and a deconcatenate coproduct. Several particular bases in the family contains the fundamental basis $\mathcal{F}_\tau$ of $\fqsym$.

Finally, in Section \ref{sec:mnrule} we prove a Murnaghan-Nakayama change of basis rule from the quasisymmetric powersum basis to the quasisymmetric fundamental basis.

Figure \ref{fig:alg} above summarizes the relationships between our new powersum bases for $\qsym$ and $\ncqsym$ and existing powersum bases in other algebras.

\section{Preliminaries}\label{sec:prelim}
Let $X=\{ x_1,x_2\ldots \}$ be a set of variables and $\alpha = (a_1,\ldots,a_k)$ be a composition of $n$. Then a generating function is quasisymmetric if, for every $k$ and $i_1<i_2<\cdots<i_k$, the coefficient of $x_{i_1}^{\alpha_1}x_{i_2}^{\alpha_2}\cdots x_{i_k}^{\alpha_k}$ is equal to the coefficient of $x_{1}^{\alpha_1}x_{2}^{\alpha_2}\cdots x_{k}^{\alpha_k}$. The set of all quasisymmetric functions is a graded Hopf algebra denoted as $\qsym$. One of the more natural bases of $\qsym$ is the \textit{quasisymmetric monomial basis} $M_\alpha$, defined as
\[
M_\alpha = \sum_{i_1<i_2<\cdots<i_k}x_{i_1}^{\alpha_1}x_{i_2}^{\alpha_2}\cdots x_{i_k}^{\alpha_k}. 
\]
 $\qsym$ is most notable for the \textit{quasisymmetric fundamental basis} $F_\alpha$ because it corresponds to the characters of the 0-Hecke algebra. To define this basis, let $\mathcal{P}$ be the refinement poset on compositions, i.e. with the cover relation $\alpha\lessdot\beta$ if $\beta = (a_1,\ldots,a_i+ a_{i+1},\ldots, a_j)$ for some $i$. Then 
\[
F_\alpha = \sum_{\beta\leq\alpha} M_\beta.
\]

Functions in $\ncqsym$ are indexed by set compositions (also called ordered set partitions), which are obtained by replacing each integer part of a composition by a set. To be more exact, a \textit{set composition} $\Phi=(B_1,\ldots,B_k)$ of $n$ is a composition of subsets of $[n]$ such that $B_i\cap B_j=\emptyset$ where $i\neq j$ and the size of $\Phi$ is $|B_1|+\ldots+|B_k|=n$. In this extended abstract we will denote a set composition $(B_1,\ldots,B_k)$ as $B_1|\ldots|B_k$, for example the set composition $(\{5\},\{1,3\},\{2\},\{4\})$ is written as $5|13|2|4$. The refinement poset $\widetilde{\mathcal{P}}$ on set compositions has the cover relation $\Phi\lessdot\Psi$ if $\Psi = (B_1,\ldots,B_i\cup B_{i+1},\ldots, B_j)$. 

Denote $\rho$ as the map from set compositions of $n$ to compositions of $n$ by $\Phi\to(|B_1|,\ldots,|B_k|)$. For example $\rho(5|13|2|4)=(1,2,1,1)$. Let $\alpha = (a_1,\ldots, a_n)$ be a composition, then $\varrho(\alpha)$ is the set composition $\Phi=(B_1,\ldots,B_k)$ where $i\in B_{|\{ a_l:a_l<a_i \}|+1}$. For example $\varrho(1,6,4,3,6)=1|4|3|25$.

Quasisymmetric function in non-commuting variables, denoted $\ncqsym$, is the space spanned by $M_\Phi$ where $M_\Phi$ is defined as 
\begin{equation*}
    M_\Phi[X]=\sum_{\alpha:\varrho(\alpha)=\Phi} x_{a_1}x_{a_2}\cdots x_{a_n}.
\end{equation*}
The map $\rho$ from set compositions to compositions induces a Hopf morphism that we will also denote as $\rho$: $\rho:\ncqsym\to \qsym$ is given by $\rho(M_\Phi)=M_{\rho(\Phi)}$.

\section{Descending Quasisymmetric Powersum Basis}\label{sec:qsym}
In this section we define a set of fillings analogous to $A(\lambda)$ in order to define the quasisymmetric powersum basis. Let $\alpha=(a_1,\ldots,a_k)$ be a composition of $n$.

\begin{defn} \label{def:strdia}
A \emph{Strict Diagonal} filling is a filling such that: 
\begin{enumerate}
    \item entry $a_{1}$ is in the upper left most corner of the matrix.
    \item entry $a_{i+1}$ is directly below $a_i$ or is in the southeast position of $a_i$ if $a_i\geq a_{i+1}$, or $a_{i+1}$ is in the southeast position of $a_i$ otherwise.
\end{enumerate}
\end{defn}

Denote the set of strict diagonal fillings with $\row(\mathsf{F})=\alpha$ as $\sd(\alpha)$. Given any filling and any integer $i$, we denote by $\sigma_i$ a permutation of the rows with entries $i$. We will only consider $\sigma_i$ such that for any two entries $a_j=a_k=i$, $\sigma_i(a_j)=a_k$ implies that $a_j$ and $a_k$ are not in the same column. Denote $\sigma=(\sigma_1,\cdots,\sigma_n)$ as a tuple of \textit{row permutations} of all the parts of $\alpha$. 
The set of all such row permutations for a filling $F$ is $\mathfrak{S}_{\mathsf{F}}$.  
Thus we define the \textit{descending quasisymmetric powersum basis} 

\begin{equation} \label{eqn:PtoMSD}
    P_\alpha =\sum_{\substack{\mathsf{F}\in \sd(\alpha)\\ \sigma \in \mathfrak{S}_{\mathsf{F}}}}M_{\col(\sigma(\mathsf{F}))}.
\end{equation}

\begin{example} \label{ex:SDfill}
Let the first two fillings below be denoted as $\mathsf{F}_1$ and $\mathsf{F}_2$, then $\sd(212)=\{\mathsf{F}_1,\mathsf{F}_2 \}$. The third and fourth fillings are $\sigma(\mathsf{F}_1)$ and $\sigma(\mathsf{F}_2)$ where $\sigma = (\text{id},(21))$. Thus $P_{(2,1,2)}=2M_{(2,1,2)}+2M_{(3,2)}$.
\[ \footnotesize
\begin{array}{|ccc}
2 & \gdot  & \gdot  \\
\gdot  & 1 & \gdot  \\
\gdot  & \gdot  & 2 \\
\hline
2 & 1 & 2
\end{array}
\qquad
\begin{array}{|ccc}
2 & \gdot  & \gdot \\
1 & \gdot   & \gdot \\
\gdot  & 2 & \gdot \\
\hline
3 & 2
\end{array}
\qquad
\begin{array}{|ccc}
\gdot  & \gdot  & 2 \\
\gdot  & 1 & \gdot  \\
2 & \gdot  & \gdot  \\
\hline
2 & 1 & 2
\end{array}
\qquad
\begin{array}{|ccc}
\gdot  & 2 & \gdot \\
1 & \gdot   & \gdot \\
2 & \gdot  & \gdot \\
\hline
3 & 2
\end{array}
\]
\end{example}

\begin{thm} \label{thm:refine}
The quasisymmetric powersum functions refine the symmetric powersum functions. In other words,
\begin{equation} \label{eqn:refine}
p_\lambda=\sum_{\alpha:\sort(\alpha)=\lambda}P_\alpha
\end{equation}
\end{thm}
 The sketch of this proof is as follows: we define a new type of filling $Q(\lambda)$ such that $p_\lambda=\sum_{\mathsf{F}\in Q(\lambda)}M_{\col(\mathsf{F})}$. Then we show that there is a bijection between $Q(\lambda) $ and the fillings of $\bigcup_\alpha \sd(\alpha)$ where $\alpha$ are all compositions such that $\sort(\alpha)=\lambda$.

\subsubsection*{Scaled Quasisymmetric Powersum Basis}

We first introduce a scaling factor so as to obtain cleaner formulas for the product and coproduct in the quasisymmetric powersum basis.

The symmetric powersums are self dual such that $\langle p_\lambda,p_\mu\rangle = z_\lambda\delta_{\lambda\mu}$ where $z_\lambda =\prod_ii^{m_i(\lambda)}m_i(\lambda)!$ where $m_i(\lambda)$ is the number of parts of size $i$ in $\lambda$. The scaled symmetric powersums are defined as $\tilde p_\lambda=\frac{1}{z_\lambda}p_\lambda$, and has the property that $\langle p_\lambda,\tilde p_\mu\rangle=\delta_{\lambda\mu}$. Analogously, we define the \textit{scaled quasisymmetric powersum basis} as $\tilde P_{\alpha} = \frac{1}{z_\alpha} P_{\alpha}$, which by Theorem \ref{thm:refine} implies that $\tilde P_\alpha$ refines $\tilde p_\lambda$. Let $Z_\alpha$ be the left Zassenhaus basis in $\nsym$ as defined in \cite{krob1997noncommutative}, then the dual of $\tilde P_\alpha$ is the left Zassenhaus basis in $\nsym$, i.e. $\langle \tilde P_\alpha,Z_\beta\rangle =\delta_{\alpha\beta}$.

\subsubsection*{Product and Coproduct}
Let $\alpha = (a_1,\ldots,a_j)$ and $\beta=(b_1,\ldots,b_k)$ be compositions of $n$ and $m$ respectively. Their concatenation is the composition $\alpha|\beta=(a_1,\ldots,a_j,b_1,\ldots,b_k)$ of $n+m$, and their shuffle, a set of compositions of $n+m$, is defined recursively by 
\begin{enumerate}
    \item $\emptyset\shuffle \alpha = \alpha \shuffle \emptyset = \alpha$
    \item $\alpha\shuffle\beta = a_1|\big((a_2,\ldots,a_j)\shuffle \beta\big) + b_1|\big( \alpha\shuffle (b_2,\ldots,b_k) \big)$.  
\end{enumerate}

\begin{thm} \label{thm:prodco}
The scaled quasisymmetric powersum basis has a shuffle product and a deconcatenate coproduct, i.e.
\begin{equation} \label{eqn:prodco}
\tilde P_\alpha \tilde P_\beta = \sum_{\gamma \in \alpha\shuffle\beta} \tilde P_\gamma,
\qquad
\Delta(\tilde P_\gamma)=\sum_{\alpha|\beta=\gamma}\tilde P_\alpha\otimes \tilde P_\beta.
\end{equation}
\end{thm}

This can be shown either using the product and coproduct of the left Zassenhaus basis, or by defining a product and coproduct on fillings.  
We give a rough sketch of the latter, as it can be generalized to the case of non-commuting variables. 

For simplicity, consider the case where all parts $c_1, \dots ,c_k$ in $\gamma$ are distinct, so $\tilde P_\gamma = \sum_{\mathsf{F}\in \sd(\gamma)}M_{\col(\mathsf{F})}$. Then 
\[ \Delta(\tilde P_\gamma)= \sum_{\substack{\mathsf{F}\in \sd(\gamma)\\ \alpha'|\beta' = \col(\mathsf{F})}}M_{\alpha'}\otimes M_{\beta'}.\] Meanwhile, \[ \sum_{\alpha|\beta=\gamma}\tilde P_\alpha\otimes \tilde P_\beta = \sum_{\substack{\mathsf{F}_1\in \sd(\alpha),\mathsf{F}_2\in \sd(\beta)\\ \alpha|\beta=\gamma}}M_{\col(\mathsf{F}_1)}\otimes M_{\col(\mathsf{F}_2)}.\] To biject between the terms on the two right hand sides above, we consider the \textit{deconcatenation of fillings}: for an $\sd$ filling $\mathsf{F}$, we write $\mathsf{F}=\mathsf{F_1}|\mathsf{F_2}$ to mean a vertical line is drawn in between two columns of $\mathsf{F}$ to get two $\sd$ fillings $\mathsf{F}_1$ and $\mathsf{F}_2$ (after removing empty rows). Note that $\col(\mathsf{F}_1)|\col(\mathsf{F}_2)=\col(\mathsf{F})$ and $\row(\mathsf{F}_1)|\row(\mathsf{F}_2)=\row(\mathsf{F})$. In \cite{duchamp2002noncommutative}, the coproduct of matrices in $\mqsym$ is defined in the exact same way (except that our fillings are transposed to their matrices).
\begin{example} The following are all four deconcatenations of one filling.
\[\footnotesize
\begin{array}{c|::cccc}
\hdashline
\hdashline
2 & 2      & \gdot  & \gdot & \gdot  \\
2 & 2      & \gdot  & \gdot & \gdot  \\
1 & \gdot  & 1      & \gdot & \gdot  \\
3 & \gdot  & \gdot  & 3     & \gdot  \\
\hline
& 4 & 1 & 3 &
\end{array}
\quad
\begin{array}{c|c::ccc}
2 & 2      & \gdot  & \gdot & \gdot  \\
2 & 2      & \gdot  & \gdot & \gdot  \\
\hdashline
\hdashline
1 & \gdot  & 1      & \gdot & \gdot  \\
3 & \gdot  & \gdot  & 3     & \gdot  \\
\hline
& 4 & 1 & 3 &
\end{array}
\quad
\begin{array}{c|cc::cc}
2 & 2      & \gdot  & \gdot & \gdot  \\
2 & 2      & \gdot  & \gdot & \gdot  \\
1 & \gdot  & 1      & \gdot & \gdot  \\
\hdashline
\hdashline
3 & \gdot  & \gdot  & 3     & \gdot  \\
\hline
& 4 & 1 & 3 &
\end{array}
\quad
\begin{array}{c|ccc::c}
2 & 2      & \gdot  & \gdot & \gdot  \\
2 & 2      & \gdot  & \gdot & \gdot  \\
1 & \gdot  & 1      & \gdot & \gdot  \\
3 & \gdot  & \gdot  & 3     & \gdot  \\
\hdashline
\hdashline
\hline
& 4 & 1 & 3 &
\end{array}
\] 
\end{example}

A similar bijective proof of the product rule relies on the \textit{quasishuffle of fillings}, defined recursively as follows.
Let $\mathsf{F} = \mathsf{F}_1|\ldots|\mathsf{F}_j$ and $\mathsf{G}= \mathsf{G}_1|\ldots|\mathsf{G}_k$, where $\mathsf{F}_i$ and $\mathsf{G}_i$ are fillings with a single column. Then set
\begin{enumerate}
    \item $\emptyset \qshuf \mathsf{F} = \mathsf{F} \qshuf \emptyset = \mathsf{F}$ where $\emptyset$ is the empty filling
    \item $\mathsf{F}\qshuf\mathsf{G} = \mathsf{F}_1|\big((\mathsf{F}_2,\ldots,\mathsf{F}_j)\qshuf \mathsf{G}\big) + \mathsf{G}_1|\big( \mathsf{F}\qshuf (\mathsf{G}_2,\ldots,\mathsf{G}_k) \big) + (\mathsf{F}_1+\mathsf{G}_1)|\big( (\mathsf{F}_2,\ldots,\mathsf{F}_j) \qshuf (\mathsf{G}_2,\ldots,\mathsf{G}_k) \big)$  
\end{enumerate}
where $\mathsf{F}_i + \mathsf{G}_{i'}$ is a filling with a single column whose entries are those in $\mathsf{F}_i$ and $\mathsf{G}_{i'}$, arranged in descending order. 
Note that, for any quasishuffle of $\mathsf{F}$ and $\mathsf{G}$, its column reading is a quaishuffle of $\col(\mathsf{F})$ and $\col(\mathsf{G})$, as in the product of monomial basis elements, and its row reading is a shuffle of $\row(\mathsf{F})$ and $\row(\mathsf{G})$, leading to a shuffle product for the scaled quasipowersum basis.

\begin{example}
The quasishuffle of two fillings.

\footnotesize
\begin{align*}
    \begin{array}{c|cc}
    2 & 2      & \gdot \\
    3 & \gdot  & 3     \\
    \hline
    & 2 & 3
    \end{array}
    \, &\qshuf \,
    \begin{array}{c|c}
    4 & 4 \\
    1 & 1 \\
    \hline
      & 5 
    \end{array} \\
    &=
    \begin{array}{c|c}
    2 & 2  \\
    \hline
      & 2 
    \end{array}
    \Bigg|\left(
    \begin{array}{c|c}
    3 & 3  \\
    \hline
      & 3 
    \end{array}
    \: \qshuf \:
    \begin{array}{c|c}
    4 & 4 \\
    1 & 1 \\
    \hline
      & 5 
    \end{array}
    \right) + \,
    \begin{array}{c|c}
    4 & 4 \\
    1 & 1 \\
    \hline
      & 5 
    \end{array}
    \Bigg| \left(
    \begin{array}{c|cc}
    2 & 2      & \gdot \\
    3 & \gdot  & 3     \\
    \hline
    & 2 & 3
    \end{array} 
    \: \qshuf \: \emptyset \right) +
    \begin{array}{c|cc}
    4 & 4 \\
    2 & 2  \\
    1 & 1  \\
    \hline
      & 7
    \end{array}
    \Bigg| \left(
    \begin{array}{c|c}
    3 & 3  \\
    \hline
      & 3 
    \end{array}
    \: \qshuf  \: \emptyset \right)      \\
    &=
    \begin{array}{c|ccc}
    2 & 2 & \gdot & \gdot \\
    3 & \gdot & 3 & \gdot  \\
    4 & \gdot & \gdot & 4 \\
    1 & \gdot & \gdot & 1 \\
    \hline
      & 2 & 3 & 5
    \end{array}
    + 
    \begin{array}{c|ccc}
    2 & 2 & \gdot & \gdot \\
    4 & \gdot & 4 & \gdot \\
    1 & \gdot & 1 & \gdot \\
    3 & \gdot & \gdot & 3 \\
    \hline
      & 2 & 5 & 3
    \end{array}
    + 
    \begin{array}{c|cc}
    2 & 2 & \gdot \\
    4 & \gdot & 4 \\
    3 & \gdot & 3  \\
    1 & \gdot & 1 \\
    \hline
      & 2 & 8
    \end{array}
    +
    \begin{array}{c|ccc}
    4 & 4 & \gdot & \gdot \\
    1 & 1 & \gdot & \gdot \\
    2 & \gdot & 2 & \gdot \\
    3 & \gdot & \gdot & 3 \\
    \hline
      & 5 & 2 & 3
    \end{array}
    +
    \begin{array}{c|cc}
    4 & 4 & \gdot \\
    2 & 2 & \gdot \\
    1 & 1 & \gdot   \\
    3 & \gdot & 3 \\
    \hline
      & 7 & 3
    \end{array}
\end{align*}
\end{example}

The product of matrices in \cite{duchamp2002noncommutative} is defined as the quasishuffle of rows; since the fillings described here are the transpose of their matrices, it is equivalent to our quasishuffle of columns. The only other difference is that the quasishuffle in this extended abstract additionally diagonalizes the fillings.

\subsubsection*{Generalizations of $P_\alpha$ Using Total Orders}

Notice that in Definition \ref{def:strdia} we put two entries in the same column if $a_i\geq a_{i+1}$. However we can generalize a strict diagonal filling according to a total order $\succeq$ where we put two entries in the same column if $a_i\succeq a_{i+1}$, then denote the analogously defined quasisymmetric powersum function as $P_\alpha^{\succeq}$. Then $P_\alpha^\succeq$ has a shuffle product, deconcatenate coproduct, and refines the symmetric powersums, as we can switch $\geq$ to $\succeq$ in relevant proofs. Note that the scaled version of $P_\alpha^{\leq}$ is the right Zassenhaus basis.

An alternative definition of $P_\alpha^{\succeq}$ uses the subposet $\mathcal{R}$ of $\mathcal{P}$, with cover relations of the form $(a_1, \ldots, a_j) \lessdot_{\mathcal{R}} (a_1,\ldots,a_i+ a_{i+1},\ldots, a_j)$ where $a_i \succeq a_{i+1}$. Then $P_\alpha^{\succeq} = \sum_{\alpha\leq_\mathcal{R}\beta}C_{\alpha\beta}M_\beta$ for some numbers $C_{\alpha\beta}$ which we won't address here.

\section{Descending Quasisymmetric Powersum Basis in Non-commuting Variables} \label{sec:ncqsym}
To generalize the above results to non-commuting variables, we consider fillings containing sets instead of integers, which we denote by $\tilde{\mathsf{F}}$. Analogous to above, the column (respectively, row) reading is a set composition recording the union of all sets in each column (respectively, row), denoted as $\col(\tilde{\mathsf{F}})$ (respectively $\row(\tilde{\mathsf{F}})$).

To define a set analog of a $\sd$ filling, we need to replace the condition $a_i \geq a_{i+1}$ in Definition \ref{def:strdia} by a comparison of blocks in the set composition. To this end, let $\min(B_i)$ be the smallest integer in $B_i$, and define $A>_{\widetilde{\mathcal{D}}}B$ if either $|A|>|B|$, or $|A|=|B|$ and $\min(A)<\min(B)$.

\begin{defn}
A \emph{Labelled Diagonal Descending (LDD) filling} of $\Phi$ is an assignment of the blocks $B_i$ to entries of a matrix such that
\begin{enumerate}
    \item $B_i$ is in row $i$ and $B_1$ is in the first column.
    \item $B_{i+1}$ can be in the same column as $B_i$ if $B_i>_{\widetilde{\mathcal{D}}}B_{i+1}$. Otherwise, $B_{i+1}$ is in the column to the right of $B_i$.
\end{enumerate}
\end{defn}
\begin{example} \label{ex:ldd}
The first two fillings are all the fillings in $\ldd(34|15|2)$, the latter four fillings are all the fillings in $\ldd(15|34|2)$.
\[ \footnotesize
\begin{array}{|ccc}
34 & \gdot  & \gdot  \\
\gdot  & 15 & \gdot  \\
\gdot  & \gdot  & 2 \\
\hline
34 |& 15| & 2
\end{array}
\qquad
\begin{array}{|ccc}
34 & \gdot  & \gdot  \\
\gdot  & 15 & \gdot  \\
\gdot  & 2  & \gdot  \\
\hline
34 |& 125| & 
\end{array}
\qquad
\begin{array}{|ccc}
15 & \gdot  & \gdot  \\
\gdot  & 34 & \gdot  \\
\gdot  & \gdot  & 2 \\
\hline
15 |& 34| & 2
\end{array}
\qquad
\begin{array}{|ccc}
15 & \gdot  & \gdot  \\
\gdot  & 34 & \gdot  \\
\gdot  & 2  & \gdot  \\
\hline
15 |& 234| & 
\end{array}
\qquad
\begin{array}{|ccc}
15 & \gdot  & \gdot  \\
34  & \gdot & \gdot  \\
\gdot  & 2  & \gdot  \\
\hline
1345 |& 2| & 
\end{array}
\begin{array}{|ccc}
15 & \gdot  & \gdot  \\
34 & \gdot & \gdot  \\
2  & \gdot  & \gdot \\
\hline
12345 &  & 
\end{array}
\]
\end{example}

\begin{defn}
The descending quasisymmetric powersum basis in NCQSym is defined as
$$P_\Phi = \sum_{\tilde{\mathsf{F}} \in \ldd(\Phi)}M_{\col(\tilde{\mathsf{F}})}. $$
\end{defn}

\begin{example}
The previous example yields $P_{34|15|2}=M_{34|15|2}+M_{34|125}$ and $P_{15|34|2}=M_{15|34|2}+M_{15|234}+M_{1345|2}+M_{12345|}$.
\end{example}

As mentioned in the Introduction, the analogs of $\sym$ and $\qsym$ in non-commuting variables are $\ncsym$ and $\ncqsym$. $\ncsym$ is indexed by set partitions $\phi$. This space has a powersum basis denoted $p_\phi$; further information about $\ncsym$ can be found in \cite{Rosas2006SymmetricVariables}. Just as in Theorem \ref{thm:refine}, quasisymmetric powersum functions in non-commuting variables refines the symmetric powersum functions in non-commuting variables. Let $\sort$ denote the map from set compositions of $n$ to set partitions of $n$ by forgetting the order of the blocks.

\begin{thm} \label{thm:toncsym}
The $P_\Phi$ basis refines the symmetric powersum function in NCSym. In other words
$$ p_\phi = \sum_{\Phi:\sort(\Phi)=\phi} P_\Phi.$$
\end{thm}

The proof of this theorem is analogous to the proof of Theorem \ref{thm:refine}.

\subsubsection*{Product and Coproduct}
Consider disjoint subsets $A_1, \dots, A_j, B_1, \dots, B_k$ of integers. Define the shifted shuffle of $\Phi=B_1|\cdots|B_k$ and $\Psi=A_1|\cdots|A_j$ recursively by
 \begin{enumerate}
    \item $\emptyset\shuffle \Phi = \Phi \shuffle \emptyset = \Phi$
    \item $\Phi\shuffle\Psi = B_1|\big((B_2,\ldots,B_j)\shuffle \Psi\big) + A_1|\big( \Phi\shuffle (A_2,\ldots,A_k) \big)$.  
\end{enumerate}
Let $\Psi\uparrow n$ denote the operation of adding $n$ to every element of $A_i$. And, if $m=|B_1|+|B_2|+\ldots+|B_k|$, let the standardization of $\Phi$, denoted $st(\Phi)$, be the set composition of $m$ where the relative order of the elements across all the $B_i$ is preserved. 
 
 \begin{thm}
The quasisymmetric powersum basis of $\ncqsym$ has a shifted shuffle product, and its coproduct is given by the standardization of the deconcatenation of blocks, i.e. if $\Phi=B_1|\cdots|B_k$ is a set composition of $n$,
\[
P_\Phi P_\Psi = \sum_{\Gamma \in \Phi\shuffle(\Psi \uparrow n)}P_\Gamma,
\qquad
\Delta(P_\Phi)=\sum_{i=0}^{k}P_{st(B_1|\cdots|B_i)}\otimes P_{st(B_{i+1}|\cdots|B_k)}.
\]
\end{thm}

The proof is analogous to the sketch above for Theorem \ref{thm:prodco}, the $\qsym$ case.

\subsubsection*{Projection onto QSym}

Let $\Phi = B_1|\cdots|B_l$ be a set composition of $n$ with $\sort(\rho(\Phi))=1^{\alpha_1} 2^{\alpha_2} \cdots j^{\alpha_j}$ and $\sigma=(\sigma_1,\dotsc,\sigma_j)$ be a permutation of $\mathfrak{S}_{\alpha_1}\times \mathfrak{S}_{\alpha_2}\cdots\times \mathfrak{S}_{\alpha_j}=\mathfrak{S}_\alpha$, where $\sigma_i$ is the permutation of blocks $B_k$ with $\len(B_k)=i$. For example, $((2,3),(1,2))(4|25|7|13|6)=4|13|6|25|7$.

\begin{thm} \label{thm:toqsym} 
For any composition $\alpha$ such that $\sort(\alpha) = 1^{\alpha_1} 2^{\alpha_2} \cdots j^{\alpha_j}$, and set composition $\Phi$ such that $\rho(\Phi)=\alpha$,
\begin{equation} \label{eq:P to P}
P_\alpha = \sum_{\sigma\in\mathfrak{S}_\alpha}\rho(P_{\sigma(\Phi)}).
\end{equation}
\end{thm}
From Example \ref{ex:ldd}, $P_{(2,2,1)}=P_{34|15|2}+P_{15|34|2}$.

\subsubsection*{Generalizations of $P_\Phi$ Using Total Orders}

With careful work we can generalize a $\ldd$ filling by replacing $>_{\tilde{\mathcal{D}}}$ by a different total ordering $\rhd$ on disjoint integer sets. The resulting powersum $P^{\rhd}_\Phi$ refines the symmetric powersums in non-commuting variables, has a shifted shuffle product, and a deconcatenate coproduct. Let $B$ and $A$ be two disjoint sets of different sizes. If there exists a total ordering on integers $\succ$, such that $B\rhd A$ implies $|B|\succ|A|$, then there exists an analog of Theorem \ref{thm:toqsym} such that sums of $P^\rhd_\Phi$ project to $P^\succ_\alpha$.

\subsubsection*{Image of FQSym}
$\fqsym$ is the Hopf algebra of permutations first introduced in \cite{malvenuto1995duality}, which is a subalgebra of $\ncqsym$ and has a fundamental basis $\mathcal{F}_\tau$. 

\begin{thm} \label{thm:fima}
The image of the $\mathcal{F}$ basis on $\ncqsym$ is the quasisymmetric powersum function in non-commuting variables indexed by a singleton set composition, i.e. for a permutation $\tau$
\begin{equation}\label{eqn:fima}
    \mathcal{F}_\tau = P_{\tau^{-1}(1|2|\cdots|n)}
\end{equation}
\end{thm}
Taking the image under $\rho$ of both sides of Equation \ref{eqn:fima} gives the equation $\rho(P_{\tau^{-1}(1|2|\cdots|n)})=F_{\des(\tau)}$, where $\des$ denotes the descent composition. It is then natural to ask for the expansion of $\rho(P_{\Phi})$ into quasisymmetric fundamental functions when the blocks of $\Phi$ are not all singletons.
To do so, we define two operations on set compositions. Let $\Phi = B_1|\ldots|B_k$ be a set composition of $n$ with $\rho(\Phi)=(a_1,\ldots,a_k)$. Break $\rho(\Phi)$ when $B_i\leq_{\widetilde{\mathcal{D}}}B_{i+1}$ so that $\rho(\Phi)=\alpha_1|\alpha_2|\cdots|\alpha_l$ where $\alpha_i=(a_{j_1},\ldots,a_{j_1+j_2})$ with  $B_{j_1-1}\leq_{\widetilde{\mathcal{D}}}B_{j_1}\geq_{\widetilde{\mathcal{D}}}B_{j_1+1}\geq_{\widetilde{\mathcal{D}}}\cdots\geq_{\widetilde{\mathcal{D}}} B_{j_1+j_2}\leq_{\widetilde{\mathcal{D}}}B_{j_1+j_2+1}  $. Then we define $\rho_C(\Phi)=(|\alpha_1|,|\alpha_2|,\cdots,|\alpha_l|)$ and $\rho_I(\Phi)= \alpha_1'|\alpha_2'|\cdots|\alpha_l'$, where $\alpha_i'$ is the transpose of $\alpha_i$.
For example, $\rho_C(2|5|14|36|7)=(2,5)$ and $\rho_I(2|5|14|36|7)=(2,1,2,2)$. 

\begin{thm} \label{thm:pwrsmim}
Let $\Phi$ be a set composition and $\mu_{\widetilde{\mathcal{P}}}$ be the M\"obius function on the refinement poset $\widetilde{\mathcal{P}}$, then
\[
\rho(P_\Phi)=\sum_{\rho_I(\Phi)\leq\alpha\leq\rho_C(\Phi)}\mu_{\widetilde{\mathcal{P}}}(\alpha,\rho_C(\Phi))F_\alpha.
\]
\end{thm}

\section{Change of basis} \label{sec:mnrule}
The Murnaghan Nakayama rule has two perspectives in $\sym$: it is the change of basis rule from symmetric powersums to Schur functions, and it is the product rule of a powersum function and a Schur function. These two perspectives are related due to the fact that symmetric powersums are multiplicative. In $\qsym$, bases are not multiplicative and so we will define the Murnaghan Nakayama rule based off of the first perspective according to \cite{Stanley2012EnumerativeCombinatorics}.
Let $\mu$ be a partition and $k$ an integer, then the product is
\begin{equation}
    s_\mu p_k = \sum_\lambda(-1)^{\hgt(\lambda/\mu)}s_\lambda.
\end{equation}
where $\lambda/\mu$ is a border strip of size $r$. This becomes a change of basis rule if we set $\mu$ to be the empty partition. Now we move onto the analog in $\qsym$.

Recall that a ribbon is a tableau such that there are no $2\times2$ boxes; for this paper we will allow ribbons to be disconnected.
 \begin{defn}
Let $\alpha$ be a composition of $n$. A filling of the ribbon of $\alpha$ is a \textit{Standard Descent Ribbon} if the Ribbon filling is standard (i.e. the numbers 1 to $n$ each appear once) and increasing from left to right, and decreasing from top to bottom.
\end{defn}
For this paper we will need to use tuples of standard descent ribbons, for example 
\[ \footnotesize
\left( \;
\ytableausetup{notabloids}
\begin{ytableau}
3 \\
\none & 1 & 4 \\
\none & \none & 2 \\
\end{ytableau}
\; , \;
\ytableausetup{notabloids}
\begin{ytableau}
2 \\
\none & 1  \\
\end{ytableau}
\; , \; 
\ytableausetup{notabloids}
\begin{ytableau}
1 & 2 \\
\end{ytableau}
\; ,\; 
\begin{ytableau}
1  \\
\end{ytableau}
\right).
\]

Let $\beta=(b_1,\ldots,b_l)$ and $\alpha=(a_1,\ldots,a_k)\leq\beta$ be compositions of $n$. $D(\beta,\alpha)=(R_1,R_2,\ldots,R_n)$ is a tuple of ribbons made from the following algorithm.

\begin{enumerate}
    \item[RT1] Let the height be equal to 0 initially and let $R$ be the ribbon of $\beta$. 
    \item[RT2] For $i$ in $[1,\ldots,k]$ one at a time, insert a box into $R_{a_i}$ in one of the following three ways:
    \begin{enumerate}
    \item[RTa1] If $a_{i-1}\neq a_i$ (or if $i=1$), then insert a box in the southeast position of the bottom right most box.
    \item[RTa2] If $a_{i-1}=a_i$ and there exists a $b_j$ such that $(b_j,\ldots,b_l)$ forms the ribbon $R$, then insert a box in the south position of the bottom right most box.
    \item[RTa3] If $a_{i-1}=a_i$ and there doesn't exist a $b_j$ such that $(b_j,\ldots,b_l)$ forms the ribbon $R$, then insert a box in the east position of the bottom right most box.
    \item[RTb] Remove the first $a_i$ boxes of $R$ and set this as $R$. The number of rows removed from $R$, minus 1, is added to the height and becomes the new height.
\end{enumerate}
\end{enumerate}

We denote the height in the algorithm as $\hgt(\beta,\alpha)$. Finally, denote by $\sdr(\alpha,\beta)$ the set of all standard descent ribbon fillings of $D(\alpha,\beta)$. We define two more operations on compositions, the analog of $\rho_C$ and $\rho_I$ of the previous section, before getting to the main theorem. Let $\alpha = (a_1,\ldots,a_k)$ be a composition of $n$, then 
break $\alpha$ at its ascents such that $\alpha=\gamma_1|\cdots|\gamma_l$ where $\gamma_i=(a_{j_1},\ldots,a_{j_1+j_2})$ where $a_{j_1-1}<a_{j_1}\geq\ldots\geq a_{j_1+j_2}<a_{j_1+j_2+1}$. Then $C(\alpha)=(|\gamma_1|,\ldots,|\gamma_l|)$ and $I(\alpha)=\gamma_1'|\cdots|\gamma_l'$.
For example, $I(1,2,1,1)=(1,1,2,1)$ and $C(1,2,1,1)=(1,4)$. 

\begin{thm}
Let $\alpha$ be a composition, then 
\begin{equation}
    P_\alpha = \sum_{I(\alpha)\leq\beta\leq C(\alpha)}(-1)^{\hgt(\beta,\alpha)}|\sdr(\beta,\alpha)|F_\beta
\end{equation}
\end{thm}

This follows from Theorem \ref{thm:toqsym} and \ref{thm:pwrsmim}. Theorem \ref{thm:toqsym} provides a way to expand a quasisymmetric powersum function to a sum of quasisymmetric powersum functions in non-commuting variables. Theorem \ref{thm:pwrsmim} takes this sum and expands it in terms of the quasisymmetric fundamental basis.

\begin{example}
Let $\alpha = (1,2,1,1)$, then to find the coefficient of $F_{(1,1,3)}$ we calculate $D(\beta,\alpha)$ when $\beta=(1,1,3)$:

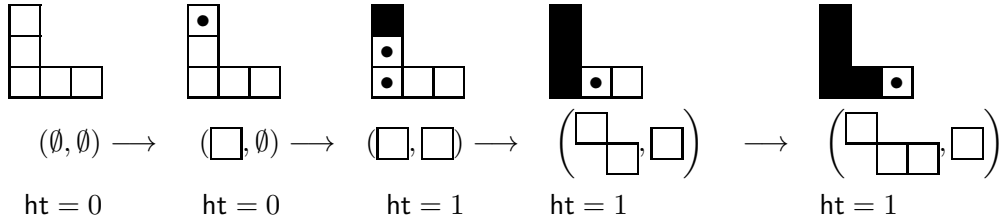
\begin{align*}
    \ydiagram{1,1,3} & & \ydiagram[\bullet]{1}*{1,1,3} & & \ydiagram[*(black)]{1}*[\bullet]{0,1,1}*{1,1,3} & & &\ydiagram[*(black)]{1,1,1}*[\bullet]{0,0,1+1}*{1,1,3}& & & &\ydiagram[*(black)]{1,1,2}*[\bullet]{0,0,2+1}&\\
    (\emptyset,\emptyset) &\longrightarrow& \left(\ydiagram{1},\emptyset \right) &\longrightarrow& \left(\ydiagram{1},\ydiagram{1} \right)&\longrightarrow& &\left(\ydiagram{1,1+1},\ydiagram{1} \right)& &\longrightarrow& &\left(\ydiagram{1,1+2},\ydiagram{1} \right)&\\
    \hgt=0 & & \hgt=0 & & \hgt=1 & & &\hgt=1& & & &\hgt=1& 
\end{align*}
and thus $\hgt(\beta,\alpha)=1$. The fillings of $D(\beta,\alpha)$ are 
\[\footnotesize
\left(
\ytableausetup{notabloids}
\begin{ytableau}
1 \\
\none & 2 & 3 \\
\end{ytableau}
\;,\;
\begin{ytableau}
1 \\
\end{ytableau}
\right)\;,\;\left(
\ytableausetup{notabloids}
\begin{ytableau}
2 \\
\none & 1 & 3 \\
\end{ytableau}
\;,\;
\begin{ytableau}
1 \\
\end{ytableau}
\right)\;,\;\left(
\ytableausetup{notabloids}
\begin{ytableau}
3 \\
\none & 1 & 2 \\
\end{ytableau}
\;,\;
\begin{ytableau}
1 \\
\end{ytableau}
\right)
\]
So the coefficient of $F_{(1,1,3)}$ is $-3$. Repeat for compositions in the interval $[I(\alpha),C(\alpha)]=[(1,1,2,1),(1,4)]$ to get $P_{(1,2,1,1)}=-3F_{(1,1,2,1)}-3F_{(1,1,3)}+3F_{(1,3,1)}+3F_{(1,4)}$.
\end{example}
\section*{Acknowledgments}
The author would like to thank his advisor Amy Pang for all the help and guidance with the development of this project and the author as a mathematician. The author would like to thank Aaron Lauve for this project, a lot of the code used and guidance. The author would also like to thank Travis Scrimshaw for all of his helpful comments. Finally the author would like to thank the Sage community.
\printbibliography

\end{document}